


\def\hexnumber#1{\ifcase#1 0\or1\or2\or3\or4\or5\or6\or7\or8\or9\or
        A\or B\or C\or D\or E\or F\fi }

\font\teneuf=eufm10
\font\seveneuf=eufm7
\font\fiveeuf=eufm5
\newfam\euffam
\textfont\euffam=\teneuf
\scriptfont\euffam=\seveneuf
\scriptscriptfont\euffam=\fiveeuf


\font\tenmsx=msam10
\font\sevenmsx=msam7
\font\fivemsx=msam5
\font\tenmsy=msbm10
\font\sevenmsy=msbm7
\font\fivemsy=msbm5
\newfam\msxfam
\newfam\msyfam
\textfont\msxfam=\tenmsx  \scriptfont\msxfam=\sevenmsx
  \scriptscriptfont\msxfam=\fivemsx
\textfont\msyfam=\tenmsy  \scriptfont\msyfam=\sevenmsy
  \scriptscriptfont\msyfam=\fivemsy
\edef\msx{\hexnumber\msxfam}

\mathchardef\upharpoonright="0\msx16
\let\restriction=\upharpoonright

\def\restrict{{\restriction}}

\def\qed{{\vcenter{\hrule height.4pt \hbox{\vrule width.4pt height5pt
 \kern5pt \vrule width.4pt} \hrule height.4pt}}}
\def\notin{{\in}\kern-5.5pt / \kern1pt}
\def\ok{\vbox{\hrule height 8pt width 8pt depth -7.4pt
    \hbox{\vrule width 0.6pt height 7.4pt \kern 7.4pt \vrule width 0.6pt height 7.4pt}
    \hrule height 0.6pt width 8pt}}
\def\nt{{\leq}\kern-1.5pt \vrule height 6.5pt width.8pt depth-0.5pt \kern 1pt}
\def\sd{{\times}\kern-2pt \vrule height 5pt width.6pt depth0pt \kern1pt}
\def\zp#1{{\hochss Y}\kern-3pt$_{#1}$\kern-1pt}

\def\sm{{\smallskip}}

\def\sub{\subseteq}

\def \o {\omega }

\font\capit=cmcsc10 scaled\magstep0

\overfullrule=0pt
\openup1.5\jot
\def\sm{\smallskip}

\def \o {\omega } \def \sub {\subseteq } 
\def\t{\tau}

\def \r {\rho } \def \k {\kappa } 
\def\g{\gamma}
\def\l{\lambda } \def\th {\theta} 
\def\d{\delta}
\def\n{\nu}

\font\sgross=cmbx10 scaled \magstep2
  \def\z{\zeta}
\def\a{\alpha } \def\G{\Gamma}
\def\b{\beta}
\def\m{\mu}

\noindent {\sgross On tightness and depth in superatomic Boolean algebras}

\bigskip 

\noindent Saharon Shelah\footnote{$^1$}{Supported by the
Basic Research Foundation of the Israel Academy of Sciences;
publication 663.}

\smallskip

\item{}{Institute of Mathematics, Hebrew University, Givat Ram, 91904
Jerusalem, ISRAEL}

{\it e-mail:} shelah@math.huji.ac.il

\smallskip

\noindent Otmar Spinas\footnote{$^2$}{Partially supported by the
Alexander von Humboldt Foundation and grant 2124-045702.95/1 of the Swiss
National Science Foundation.}

\smallskip

Mathematik, ETH-Zentrum, 8092 Z\"urich, SWITZERLAND

{\it e-mail:} spinas@math.ethz.ch

\bigskip 

{\narrower

{\noindent ABSTRACT: We introduce a large cardinal property which is
consistent with $L$ and show that for every superatomic Boolean
algebra $B$ and every cardinal $\lambda$ with the large cardinal
property, if tightness$^+(B)\geq \lambda ^+$ then depth$(B)\geq
\lambda$. This improves a theorem of Dow and Monk.

}}

\bigskip

\bigskip

In [DM, Theorem C], Dow and Monk have shown that if $\lambda$ is a
Ramsey cardinal (see [J, p.328]) then every superatomic Boolean
algebra with tightness at least $\lambda ^+$ has depth at least
$\lambda$. Recall that a Boolean algebra $B$ is {\it superatomic} iff
every homomorphic image of $B$ is atomic. The {\it depth} of $B$ is
the supremum of all cardinals $\lambda$ such that there is a sequence
$(b_\alpha :\alpha <\lambda )$ in $B$ with $b_\beta <b_\alpha $ for
all $\alpha <\beta <\lambda$ (a {\it well-ordered chain} of length
$\lambda $). Then depth$^+$ of $B$ is the first cardinal $\lambda$
such that there is no well-ordered chain of length $\lambda$ in
$B$. The {\it tightness} of $B$ is the supremum of all cardinals
$\lambda$ such that $B$ has a {\it free} sequence of length $\lambda$,
where a sequence $(b_\alpha :\alpha <\lambda )$ is called {\it free}
provided that if $\Gamma$ and $\Delta$ are finite subsets of $\lambda$
such that $\alpha <\beta$ for all $\alpha \in \Gamma$ and $\beta\in
\Delta$, then $$\bigcap_{\alpha \in \Gamma}-b_\alpha \cap
\bigcap_{\beta\in\Delta} b_\beta \ne 0.$$ By tightness$^+(B)$ we denote the
first cardinal $\lambda$ for which there is no free sequence of length
$\lambda$ in $B$. 

For $b\in B$ we sometimes write $b^0$ for $-b$ and $b^1$ for $b$.

We improve Theorem C from [DM] in two directions. We introduce a large
cardinal property which is much weaker than Ramseyness and even
consistent with $L$ (the constructible universe) and show that in
Theorem C from [DM] it suffices to assume that $\lambda$ has this
property. Moreover we show that it suffices to assume
tightness$^+(B)\geq \lambda ^+$ instead of tightness$(B)\geq \lambda
^+$ to conclude that depth$(B)\geq \lambda$. In particular we get:

\sm

{\bf Theorem 1.} {\it Suppose that $0^\sharp$ exists. Let $B$ be a
superatomic Boolean algebra in the constructible universe $L$, and let
$\lambda$ be an uncountable cardinal in $V$. Then in $L$ it is true
that tightness$^+(B)\geq \lambda ^+$ implies that depth$^+(B)\geq
\lambda$.}

\sm

For the theory of $0^\sharp$ see [J, $\S30$]. Note that $\lambda$ as
in Theorem 1 is a limit cardinal in $L$, hence it suffices to show
that in $L$, depth$(B)\geq \kappa$ for all cardinals $\kappa
<\lambda$. As was the case with the proof of Theorem C of [DM], we
can't show that under the assumptions of Theorem 1, depth$(B)=\lambda$
is attained, i.e. that there is a well-ordered chain of length
$\lambda$.

\sm

For the proof we consider the following large cardinal property:

\sm

{\bf Definition 2.}  Let $\l,\,\k,\,\th$ be infinite cardinals, and
let $\g$ be an ordinal. The relation $R_\g (\l,\k,\th )$ is defined as
follows:

{\narrower{\noindent For every $c:[\l]^{<\o}\rightarrow \th$ there
  exists $A\sub \l$ of order-type $\g$, such that for every $u\in
  [A]^{<\o}$ there exists $B\sub\l$ of order-type $\k$ such that
  $\forall w\in [B]^{|u|}\quad c(w)=c(u).$

}}

\sm

{\bf Lemma 3.} {\it Assume $R_\g(\l,\k,\th)$, where $\g$ is a limit
  ordinal. For every $c:[\l]^{<\o}\rightarrow \th$ there exists
  $A\sub\l$ as in the definition of $R_\g(\l,\k,\th)$ such that
  additionally $c\restrict [A]^n$ is constant for every $n<\o$.}

\sm

{\it Proof:} Define $c'$ on $[\l]^{<\o}$ by $$c' \{ \beta _0,\dots
,\beta _{n-1} \} = \{ (v,c\{ \beta _i:i\in v \} ) :v\sub n\}.$$ 
As $\th$ is infinite we
can easily code the values of $c'$ as ordinals in $\th$ and therefore
apply $R_\g(\l,\k,\th)$ to it. We get $A\sub\l$ of order-type $\g$. We
shall prove that $c\restrict [A]^n$ is constant, for every $n<\o$. Fix
$w_1,w_2\in [A]^n$. Since $\g$ is a limit, without loss of generality
we may assume that $\max(w_1)<\min(w_2)$. Let $w=w_1\cup w_2$. By
Definition 2 there exists $B\sub\l$, o.t.$B=\k$, such that
$c'\restrict [B]^{2n}$ is constant with value $c'(w)$. Let
$(\b_\nu:\nu<\k)$ be the increasing enumeration of $B$. We have
$$c'\{\b_0,\dots,\b_{2n-1}\}=c'\{\b_n,\dots,\b_{3n-1}\}.$$ By the
definition of $c'$ we get
$$c\{\b_0,\dots,\b_{n-1}\}=c\{\b_n,\dots,\b_{2n-1}\}=:c_0.$$ This
information is coded in $c'\{\b_0,\dots,\b_{2n-1}\}$,
i.e. $$(\{0,\dots,n-1\},c_0),\, (\{n,\dots,2n-1\},c_0)\in
c'\{\b_0,\dots,\b_{2n-1} \}.$$ As $c'\{\b_0,\dots,\b_{2n-1}\}=c'(w)$
we conclude $c(w_1)=c(w_2)=c_0$. \hfill {$\qed$}

\sm

{\bf Theorem 4.} {\it Assume $R_\g(\l,\k,\o)$, where $\g$ is a limit
  ordinal. If $B$ is a Boolean algebra and $(a_\nu:\nu<\l)$ is a
  sequence in $B$, then one of the following holds: 

\item{(a)} there exists $A\sub\l$, o.t.$(A)=\g$, such that
  $(a_\n:\n\in A)$ is independent;

\item{(b)} there exist $n<\o$ and strictly increasing sequence
  $(\b_\n:\n<\kappa)$ in $\l$ such that, letting
  $$b_\n=\bigcup_{k<n}\bigcap_{l<n} a_{\b_{n^2\n+nk+l}},\eqno(\ast)$$
    we have that $(b_\n:\n<\k)$ is constant;

\item{(c)} there exists a strictly decreasing sequence in $B$ of
  length $\k$.}

\sm

{\bf Corollary 5.} {\it Assume $R_\g(\l,\k,\o)$, where $\g$ is a limit
  ordinal. If $B$ is a superatomic Boolean algebra, then
  tightness$^+(B)>\l$ implies Depth$^+(B)> \k$. }

\sm

{\it Proof of Corollary 5:} Let $(a_\n:\n<\l)$ be a free sequence in
$B$. As a superatomic Boolean algebra does not have an infinite
independent subset, (a) is impossible. Suppose (b) were true. Define
$b_\n$ as in $(\ast)$. Clearly we have $$-b_\n\geq
\bigcap_{k,l<n}a^0_{\b_{n^2\n+nk+l}},\hbox{ and}$$
$$b_\n\geq \bigcap_{k,l<n}a_{\b_{n^2\n+nk+l}}.$$ Hence if $\n<\mu$ and
$b_\n=b_\mu$ we obtain $$0=-b_\n\cap b_\mu\geq
\bigcap_{k,l<n}a^0_{\b_{n^2\n+nk+l}} \cap
\bigcap_{k,l<n}a_{\b_{n^2\mu +nk+l}}.$$ This contradicts freeness of
$(a_\n:\n<\k)$. We conclude that (c) must hold. \hfill {$\qed$}

\sm

{\it Proof of Theorem 4:} Define $c:[\l]^{<\o}\rightarrow [{^{<\omega
    }2}]^{<\o}$ by $$c\{\b_0<\dots <\b_{n-1}\}=\{\eta\in
    {^n2}:\bigcap_{i<n}a_{\b_i}^{\eta(i)} =0\}.$$ Note that
    $c\{\b_0<\dots <\b_{n-1}\}=c\{\a_0<\dots <\a_{n-1}\} $ implies
    that $\{a_{\b_0},\dots,a_{\b_{n-1}}\}$ and
    $\{a_{\a_0},\dots,a_{\a_{n-1}}\}$ have the same quantifier-free
    diagram, i.e. for every quantifier-free formula
    $\phi(x_0,\dots,x_{n-1})$ in the language of Boolean algebra,
    $$B\models\phi[a_{\b_0},\dots,a_{\b_{n-1}}] \Leftrightarrow
    B\models\phi[a_{\a_0},\dots,a_{\a_{n-1}}] .$$ Let $A\sub \l$ be as
    guaranteed for $c$ by $R_\g(\l,\k,\o)$. By Lemma 3 we may assume
    that $c\restrict [A]^n$ is constant, for every $n<\o$.

If $(a_\a:\a\in A)$ is independent, we are done. Therefore we may
assume that this is false. For $m<\o$ define $$\G_m=\{\eta\in
{^m2}:\exists \{\b_0<\dots <\b_{m-1}\}\sub A\quad
\bigcap_{i<m}a^{\eta(i)}_{\b_i} =0\}.$$ By assumption, in the
definition of $\G_m$ the existential
quantifier can be replace by a universal one to give the same
set. There exists $m<\o$ such that $\G_m\ne\emptyset$. Define
$$\G_m'=\{\eta\in\G_m:\hbox{ no proper subsequence of } \eta \hbox{
belongs to }\bigcup_{k<m}\G_k\}.$$ By Kruscal's Theorem [K], we have
that $\bigcup_{m<\o}\G_m'$ is finite. Let $n^*$ be minimal such that
$\bigcup_{m<\o}\G_m'=\bigcup_{m<{n^*}}\G_m'$. Then clearly we have
that for every $m<\o$ and $\eta\in\G_m$, $\eta$ has a subsequence in
$\bigcup_{k<n^*}\G_k'$. Let $m^*=(n^*)^2$, and let $$\t(x_0,\dots
,x_{m^*-1})=\bigcup_{l<n^*}\bigcap_{k<n^*}x_{n^*l+k}.$$

\sm

{\bf Claim 1.}  {\it If $\eta\in {^{m^*}2}$, $t\in\{ 0,1\}$, and in the
Boolean algebra $\{0,1\}$, $\t[\eta(0),\dots,\eta(m^*-1)]=t$, then
$|\{i<m^*:\eta(i)=t\}|\geq n^*$.}\hfill {$\qed$}

\sm

Let $(\beta_\n : \n <\gamma)$ be the strictly increasing enumeration
of $A$, and define $$b_\n=\t[a_{\b _{m^*\n}},a_{\b_{m^*\n
+1}},\dots,a_{\b_{m^*\n +m^*-1}}],$$ for every $\n<\g$, where the
evaluation of $\t$ takes place in $B$, of course. It is easy to see
that the sequence $(b_\n:\n <\g)$ inherites from $(a_{\b_\n}:\n <\g)$ the
property, that any two finite subsequences of same length have the
same quantifier-free diagram.

\sm

{\bf Claim 2.} {\it If $\eta\in \G_n$, then
$\bigcap_{i<n}b^{\eta(i)}_i=0$.}

\sm

{\it Proof of Claim 2:} Otherwise there exists an ultrafilter $D$ on
$B$ such that $\bigcap_{i<n}b^{\eta(i)}_i\in D$. Define $\zeta\in
{^{nm^*}2}$ by $\zeta (i)=1$ iff $a_{\b_i}\in D$. Then
$\bigcap_{i<nm^*}a_{\b_i}^{\z(i)}\in D$, and hence
$\z\not\in\G_{nm^*}$. Let $h:B\rightarrow B/D=\{ 0,1\}$ be the
canonical homomorphism induced by $D$. We calculate 
\sm

{\narrower{\noindent
$1=h(\bigcap_{i<n}b^{\eta(i)}_i)=\bigcap_{i<n}h(b_i)^{\eta(i)}
=\bigcap_{i<n}\t [h(a_{\b_{m^*i}}), \dots,h(a_{\b_{m^*(i+1)-1}})]^{\eta(i)}$
\sm
\hfill $=\bigcap_{i<n}\t[\z(m^*i),\dots,\z(m^*i +k), \dots,
\z(m^*(i+1)-1)]^{\eta(i)}.$

}} \sm We conclude that $\t[\z(m^*i),\dots,\z(m^*i+k), \dots,
\z(m^*(i+1)-1)]=\eta(i)$, for all $i<n$, and hence by Claim 1 we can
choose $j_i\in [m^*i,m^*(i+1))$ such that $\z(j_i)=\eta(i)$. Clearly
$i_0<i_1$ implies that $j_{i_0}<j_{i_1}$. But this implies
$\z\in\G_{nm^*}$, a contradiction. \hfill {$\qed_{Claim\; 2}$}

\sm

{\bf Claim 3.} {\it If $t<\o$, $\eta\in \G_n$, $0=k_0<k_1<\dots
<k_t=n$, and $\eta\restrict [k_i,k_{i+1})$ is constant for all $i<t$,
and if $\r\in {^t2}$ is defined by $\r(i)=\eta(k_i)$, then
$\bigcap_{i<t}b_i^{\r(i)}=0$.}

\sm

{\it Proof of Claim 3:} Wlog we may assume that $\eta\in \G_n'$ for
some $n<n^*$. Indeed, otherwise we can find $m<n^*$, $\eta'\in\G_m'$
and some increasing $h:m\rightarrow n$ such that
$\eta'(i)=\eta(h(i))$, for all $i<m$. Then
$\{h^{-1}[k_i,k_{i+1}):i<t\}$ equals $\{ [l_i,l_{i+1}):i<s\}$ for some
$l_0=0<l_1<\dots <l_{s-1}=m$. Note that $\eta'\restrict [l_i,l_{i+1})$
is constant, and letting $\r'\in {^s2}$ be defined by
$\r'(i)=\eta'(l_i)$, we have $\r'(i)=\r(h(i))$. Hence
$\bigcap_{i<s}b_i^{\r'(i)}=0$ implies
$\bigcap_{i<t}b_i^{\r(i)}=0$. 

Therefore we assume $\eta\in \G_n'$, for
some $n<n^*$. Suppose we had $\bigcap_{i<t}b_i^{\r(i)}>0$. Let $D$ be
an ultrafilter on $B$ containing $\bigcap_{i<t}b_i^{\r(i)}$. Let
$h:B\rightarrow B/D$ be the canonical homomorphism. Define $\z\in
{^{tm^*}2}$ such that $\z(i)=1$ iff $a_i\in D$. Hence
$\z\not\in\G_{tm^*}$. We get
$$h(\bigcap_{i<t}b_i^{\r(i)})=\bigcap_{i<t}\t[\z(im^*) ,\dots
,\z((i+1)m^*-1)]^{\r(i)} =1.$$ Hence by Claim 1, $$\forall i<t\exists
a_i\in [\{im^*,\dots ,(i+1)m^*-1\}]^{n^*}\forall j\in a_i\quad
\z(j)=\r(i).$$ Define $\m\in {^{tn^*}2}$ by $\m (j)=\r(i)$ iff $j\in
[in^*,(i+1)n^*) $. Then $\m$ is a subsequence of $\z$ and therefore
$\m\not\in\G_{tn^*}$. But also $\eta$ is a subsequence of $\m$, and
hence $\eta\not\in \G_n$, a contradiction. 

\hfill {$\qed_{Claim\; 3}$}

\sm

{\bf Claim 4.} {\it Suppose $\r\in {^t2}$ and $\bigcap_{i<t}b_i^{\r
(i)}=0$. Let $\z\in {^{m^*t}2}$ be defined such that $\z(m^*i)=\r(i)$
and $\z\restrict [m^*i,m^*(i+1))$ is constant for every $i<t$. Then
$\z\in\G_{m^*t}$. }

\sm

{\it Proof of Claim 4:} Otherwise,
$\bigcap_{i<m^*t}a_i^{\z(i)}>0$. Let $D$ be an ultrafilter containing 
$\bigcap_{i<m^*t}a_i^{\z(i)}$. Let $h:B\rightarrow B/D$ be the
canonical homomorphism. We have
$$h(\bigcap_{i<t}b_i^{\r(i)})=\bigcap_{i<t}\t[\z (m^*i),\dots
,\z(m^*(i+1)-1)]^{\r(i)} =\bigcap_{i<t}\t[\r(i),\dots
,\r(i)]^{\r(i)}=1.$$ This is a contradiction. \hfill{$\qed _{Claim\;
4}$}

\sm

Since we assume that $(a_\a:\a\in A)$ is not independent, by Claim 2
we can find $k^*<\o$ minimal such that for some $\r^*\in {^{k^*}2}$,
$\bigcap_{i<k^*} b^{\r^*(i)}_i=0$. Note that $\r^*(i+1)\ne\r^*(i)$ for
every $i<k^*-1$. Indeed, otherwise let $\z\in {^{m^*k^*}2}$ be defined
as in Claim 4. So $\z\in\G_{m^*k^*}$. By Claim 3 we can find $\r'$ of
shorter length than $\r^*$ such that
$\bigcap_{i<|\r'|}b_i^{\r'(i)}=0$, contradicting the minimal choice of
$k^*$.

Suppose first that $k^*=1$.  We conclude that $(b_\n :\n <\gamma)$
either is constantly 1 or 0. The main part of the definition of
$R_\gamma (\lambda , \kappa , \omega )$ then gives a sequence of
length $\kappa$ as desired in (b) of Theorem 4.

Secondly suppose $k^*>1$. If $\bigcap_{i<k^*-2}b_i^{\r^*(i)}\cap
b_{k^*-2}\cap b_{k^*-1}^0 =0 $ and $\bigcap_{i<k^*-2}b_i^{\r^*(i)}\cap
b_{k^*-2}^0\cap b_{k^*-1} =0 $, then $\bigcap_{i<k^*-2}b_i^{\r^*(i)}\cap
b_{k^*-2}= \bigcap_{i<k^*-2}b_i^{\r^*(i)}\cap b_{k^*-1} $, and an application
of the main part of the definition of
$R_\gamma (\lambda , \kappa , \omega )$ gives a sequence as desired in (b). 

Otherwise, if $\r^*(k^*-2)=1$ and $\r^*(k^*-1)=0$, then
$$\bigcap_{i<k^*-2}b_i^{\r^*(i)}\cap b_{k^*-2}<
\bigcap_{i<k^*-2}b_i^{\r^*(i)}\cap b_{k^*-1} $$, and applying the definition
gives (c). Similarly if $\r^*(k^*-2)=0$ and $\r^*(k^*-1)=1$. \hfill {$\qed$}

\sm

{\bf Theorem 6.} {\it Assume the following:

\item{(1)} $0^\sharp $ exists,

\item{(2)} $V\models \l$ is an uncountable cardinal,

\item{(3)} $\k,\theta<\l$, and $L\models \k$ is a regular
cardinal. 

Then $L\models R_\o(\l,\k,\theta)$.}

\sm

{\it Proof:}  Let $c:[\l]^{<\o}\rightarrow \theta$, $c\in L$, be
arbitrary. 

Let $Y$ be the set of all $w\in [\l ]^{<\o}$ such that
for every $n\leq |w|$ and $u\in [w]^n$ there exists $B\sub \l$ of
order-type $\k$ in $L$ such that  $\forall v\in [B]^n \quad
c(u)=c(v)$. Clearly $Y\in L$.

\sm

{\bf Claim 1.} {\it If in $V$ there exists $A\in [\l]^\o$ with
$[A]^{<\o}\sub Y$, then $L\models R_\o(\l,\k,\theta)$.}

\sm

{\it Proof of Claim 1:}  Let $T$ be the set of all one-to-one
sequences $\r\in {^{<\o}\l}$ with ran$(\r)\in Y$, ordered by
extension. Then $T$ is a tree and by assumption, $T$ has an
$\o$-branch in $V$. By absoluteness, $T$ has an $\o$-branch  $b$ in
$L$. Then ran$(b)$ (or some subset) witnesses $L\models
R_\o(\l,\k,\theta)$. \hfill {$\qed _{Claim\;1}$}

\sm

Let $(i_\n :\n<\l^+)$ be the increasing enumeration of the club of
indiscernibles of $L_{\l^+}$. Then $(i_\n:\n<\l)$ is the club of
indiscernibles of $L_\l$. As $c\in L_{\l^+}$ there exist ordinals
$\xi_0<\dots <\xi_{p-1}<\l\leq\xi_p<\dots <\xi_{q-1}<\l^+$ and a
Skolem term $t_c$ such that $$L_{\l^+}\models c=t_c[i_{\xi_0},\dots
,i_{\xi_{q-1}}].$$ By indiscernibility and remarkability (see [J,
p.345]) it easily follows that if $\a^*=\max\{ \xi_{p-1},\theta\}+1$,
then $c\restrict [\{i_\n :\a^*\leq \n<\l\}]^n$ is constant for every
$n<\o$, say with value $c_n$.  Let $n<\o$ be arbitrary. Let $\d_0
=i_{\a^* +\k},\; \d_1=i_{\a^*+\k+1},\dots , \d_{n-1}=i_{\a^*+\k+n-1}.$

\sm

{\bf Claim 2.} {\it For every $\a<\d_0$ there exists a limit $\d$,
$\a<\d<\d_0$, such that for all $\b_0<\dots <\b_{n-2}<\d$ the
following hold:

\item{$(\ast)_0$} $c\{ \d,\d_1,\dots , \d_{n-1}\}=c\{ \d_0,\dots
,\d_{n-1}\} (= c_n),$

\item{$(\ast)_1$} $c\{\b_0, \d,\d_2,\dots , \d_{n-1}\}=c\{ \b_0, \d_1,\dots
,\d_{n-1}\},$

\item{$(\ast)_2$} $c\{\b_0 ,\b_1, \d,\d_3,\dots ,\d_{n-1}\}= c\{\b_0
,\b_1, \d _2,\dots ,\d_{n-1}\},$

\sm

\item{} $\dots $

\sm

\item{$(\ast)_{n-1}$} $c\{ \b_0,\dots ,\b_{n-2},\d\} = c\{ \b_0,\dots
,\b_{n-2}, \d_{n-1}\}.$ }

\smallskip

{\it Proof of Claim 2:}  Let $\a<\d_0$ be arbitrary. Choose $\g<\k$
such that $\g$ is a limit and $i_{\a^*+\g}>\a$, and let $\d
=i_{\a^*+\g}$.

Then clearly $(\ast )_0$ holds. 

In order to prove $(\ast)_1$, let $\b<\d$ be arbitrary. There exist
ordinals $\n_0<\dots <\n_{k-1}<\a^*+\g$ and a Skolem term $t_\b$ such
that $$t_\b^{L_\l}[i_{\n_0},\dots ,i_{\n_{k-1}}]=\b.$$ Moreover there
exist ordinals $\m_0<\dots <\m_{l-1}<\a^*$ and a Skolem term $t$ such that 
$$L_{\l^+}\models t[i_{\mu_0},\dots
,i_{\mu_{l-1}}]=t_c[i_{\xi_0},\dots ,i_{\xi_{q-1}}]\{
t_\b[i_{\nu_0},\dots ,i_{\n_{k-1}}],\d_1,\dots
,\d_{n-1}\}.\eqno{(+)}$$ Note that all indices of occurring
indiscernibles, except for $\d_1,\dots ,\d_{n-1}$, either are at least
$\l$ or else below $\a^*+\g$. We conclude that in $(+)$, $\d_1$ can be
replaced by $\d$. The resulting statement is $$c\{ \b,\d_1,\dots
,\d_{n-1}\} = c\{ \b,\d,\d_2, \dots ,\d_{n-1}\},$$ as desired.

The proof of $(\ast)_2${---}$(\ast)_{n-1}$ is
similar. \hfill {$\qed_{Claim\; 2}$}

\sm

It is clear that the statement of Claim 2 is absolute.  Hence it is
also true in $L$. Using this we shall prove that $[\{ i_\nu :\alpha ^*
\leq \nu
<\lambda \}]^{<\o}\sub Y$. By Claim 1, this will suffice. We only have
to prove that for every $n<\o$ there exists $B\sub \l$ of order-type
$\k$ such that $B\in L$ and $\forall v\in [B]^n\quad c(v)=c_n.$ Fix
$n<\o$. Working in $L$, we construct $B$ inductively as $\{\g_\n:\n
<\k\}$.

Fix $\d_0<\d_1<\dots <\d_{n-2}<\l$ as above. Apply Claim 2 in $L$ with
$\a=0$ and obtain $\g_0\in (0,\d_0)$. Suppose we have gotten $(\g_\n
:\n<\m)$ for some $\m<\k$. Let $\g^*=\sup_{\n<\m}\g_\n +1$. Since
cf$^L(\d_0)\geq \k$ and $(\g_\n:\n<\m)\in L$, we have that
$\g^*<\d_0$. Apply Claim 2 with $\a=\g^*$ and get $\g_\m\in
(\g^*,\d_0)$.

We claim that $(\g_\n:\n<\k)$ is as desired. Indeed, let $\{
\g_{\n_0}<\g_{\n_1}< \dots <\g_{\n_{n-1}}\}$ be arbitrary. We have
\sm
$c\{\g_{\n_0},\dots ,\g_{\n_{n-1}}\} =^{(\ast)_{n-1}} c\{\g_{\n_0},\dots
,\g_{\n_{n-2}}, \d_{n-1}\}$
\sm
$=^{(\ast)_{n-2}} c\{ \g_{\n_0},\dots ,\g_{\n_{n-3}}, \d_{n-2},
\d_{n-1}\}$

\sm

$ = \dots$

\sm

$=^{(\ast)_1} c\{ \g_{\n_0},\d_1,\dots ,\d_{n-1}\}$
\sm
$=^{(\ast )_0} c_n$. \hfill {$\qed_{Theorem\;6}$}

\bigskip

\bigskip

\centerline{{\bf References}}

\smallskip

\item{[DM]} {\capit A. Dow and D. Monk,} {\it Depth, $\pi$-character,
and tightness in superatomic Boolean algebras,} Top. and its Appl. {\bf
15}(1997), 183-199.

\item{[J]} {\capit T. Jech,} {\bf Set Theory,} Academic Press, New
York, 1978.

\item{[K]} {\capit J. Kruskal,} {\it Well-quasi ordering, the tree
theorem and Vazsonyi's conjecture,} Trans. Am. Math. Soc. {\bf
95}(1960), 210-225.

\bye